\documentclass[12pt]{amsart}

\def\@typesizes{%
       \or{5}{6.5}\or{6}{7.5}\or{7}{8.5}\or{8}{11}\or{9}{12}%
       \or{10}{13}
       \or{\@xipt}{14}\or{\@xiipt}{15}\or{\@xivpt}{18}%
       \or{\@xviipt}{20}\or{\@xxpt}{24}}

\usepackage{amsthm}
\usepackage{amsmath}
\usepackage{amssymb}
\usepackage{enumerate}
\usepackage{color}
\allowdisplaybreaks
\numberwithin{equation}{section}
\numberwithin{figure}{section}

\theoremstyle{plain}
\newtheorem{theorem}{ Theorem}[section]
\newtheorem{proposition}[theorem]{ Proposition}
\newtheorem{lemma}[theorem]{ Lemma}
\newtheorem{corollary}[theorem]{ Corollary}
\newtheorem{example}[theorem]{ Example}
\newtheorem{remark}[theorem]{ Remark}
\newtheorem{definition}[theorem]{ Definition}
\newtheorem{conjecture}{ Conjecture}

\def\BET{\begin{theorem}}
\def\ENT{\end{theorem}}
\def\BEP{\begin{proposition}}
\def\ENP{\end{proposition}}
\def\BEL{\begin{lemma}}
\def\ENL{\end{lemma}}
\def\BEC{\begin{corollary}}
\def\ENC{\end{corollary}}
\def\BEE{\begin{example} \rm}
\def\ENE{\end{example}}
\def\BER{\begin{remark} \rm}
\def\ENR{\end{remark}}
\def\BED{\begin{definition} \rm}
\def\END{\end{definition}}
\def\BECJ{\begin{conjecture}}
\def\ENCJ{\end{conjecture}}

\def\bea{\begin{eqnarray}}
\def\eea{\end{eqnarray}}

\def\beas{\begin{eqnarray*}}
\def\eeas{\end{eqnarray*}}

\def\beq{\begin{equation}}
\def\eeq{\end{equation}}

\def\beal{\begin{align*}}

\def\eeal{ \end{align*} }

\def\roweq{\nonumber \\ &=& }
\def\rowleq{\nonumber \\  & \leq & }
\def\rowgeq{\nonumber \\ & \geq & }

\def\bbC{{\mathbb C}}
\def\bbD{{\mathbb D}}

\def\bbN{{\mathbb N}}

\begin{document}

\title{Solid hulls and cores of weighted $H^\infty$-spaces.}

\author{Jos\'e Bonet}
\author{Wolfgang Lusky}
\author{Jari Taskinen}


\begin{abstract}
We determine the solid hull and solid core of weighted Banach spaces $H_v^\infty$ of analytic functions 
functions $f$ such that $v|f|$ is bounded, both in the case of the holomorphic functions on the disc and on the whole complex plane, for a very general class of radial weights $v$. Precise results are presented for concrete weights on the disc that could not be treated before. It is also shown that if $H_v^\infty$ is solid, then the monomials are an (unconditional) basis of the closure of the polynomials in $H_v^\infty$. As a consequence $H_v^\infty$ does not coincide with its solid hull and core in the case of the disc. An example shows that this does not hold for weighted spaces of entire functions.
\end{abstract}

\maketitle


\section{Introduction and preliminaries}
\label{sec1}

The solid hulls of weighted $H^\infty$-type Banach spaces $H_v^\infty$ of analytic functions
on the open
unit disc $\bbD = \{ z \in  \bbC \, :\, |z| <1\} $ 
were characterized in
\cite{BT2}  for a large class of weight functions $v$, and a similar study for entire functions
was made in \cite{BT1}. In Theorem \ref{th4.1} we extend the results of \cite{BT2} by means of the calculations of certain numerical constants, which yields many novel, concrete examples of
solid hulls. At the same time, we also describe in Theorem \ref{th2.2} the solid cores of these Banach spaces in a similar way as we did in Theorem \ref{th2.1} (Theorem 2.1 in \cite{BT1}) for the solid hull.
Moreover, we prove that in the case of analytic functions on the disc the Banach space $H_v^\infty$ is always different from both its solid hull and core; see Corollary \ref{cor5.3}. However, an example is given in the case of entire functions to show that $H_v^\infty$ may coincide with both its solid hull and core, in particular it is a solid space. We also prove in Theorem \ref{th5.2} that if $H_v^\infty$ coincides with its solid hull, then the monomials are an (unconditional) basis of the closure of the polynomials in $H_v^\infty$.

To describe the results in detail, let us introduce some notation and terminology. We set $R= 1$ (for the case of holomorphic functions on the unit disc) and $R= +\infty$ (for the case of entire functions).  A {\it weight} $v$   is a continuous function  $v: [0, R[ \to ]0,  \infty [$, which is non-increasing on $[0,R[$ and satisfies $\lim_{r \rightarrow R} r^n v(r)=0$ for each $n \in \bbN$. We extend $v$ to $\bbD$ if $R=1$ and to $\bbC$ if $R=+\infty$ by $v(z):= v(|z|)$. For
such a weight $v$, we study the Banach space  $H_v^\infty$ of analytic
functions $f$ on the disc $\bbD$ (if $R=1$) or on the whole complex plane $\bbC$ (if $R=+\infty$) such that $\Vert f \Vert_v:= \sup_{|z| <R} v(z)|f(z)| <\infty$. For an analytic
function $f \in H(\{z \in \bbC ; |z| <R \})$ and $r <R$, we denote $M(f,r):= \max\{|f(z)| \ ; \ |z|=r\}$. Using the notation $O$ and $o$ of Landau, $f \in H_v^\infty$ if and only if $M(f,r)=O(1/v(r)), r \rightarrow R$. It is known that the closure of the polynomials in $H_v^\infty$ coincides with the Banach space $H_v^0$ of all those analytic functions on $\{z \in \bbC ; |z| <R \}$ such that $M(f,r)=o(1/v(r)), r \rightarrow R$. see e.g.\ \cite{BBG}. It will be clear from the context in the rest of the article when we refer to analytic functions on the disc or entire functions. Anyway, if it is necessary to distinguish at some point, we will use the notations $H_v^\infty(\bbD)$ and $H_v^\infty(\bbC)$.

We shall identify an analytic function $f(z)=\sum_{n=0}^\infty a_n z^n$  with the
sequence of its Taylor coefficients $(a_n)_{n=0}^\infty$. Let $A$, $B$ and $H$ be vector
spaces of complex sequences containing the space of all the sequences with finitely many non-zero coordinates.
The space $A$ be is \textit{solid} if $a=(a_n) \in A$ and $|b_n| \leq |a_n|$ for each $n$ implies
$b=(b_n) \in A$. The \textit{solid hull of $A$} is
$$S(A):= \{ (c_n) \, : \, \exists (a_n) \in A \ \mbox{such that} \ |c_n| \leq |a_n| \ \forall n \in \bbN \}. $$
It coincides with the smallest solid space containing $A$.

\noindent The \textit{solid core of $A$} is

$$s(A):=\{ (c_n) \, : \, (c_na_n) \in A \ \forall (a_n) \in \ell_\infty \}.$$

\noindent The \textit{set of multipliers form $A$ into $B$} is

$$(A,B):= \{ c=(c_n) \, : \, (c_na_n) \in B \ \forall (a_n) \in A \}. $$

The following facts are also known, see \cite{AS}:  1. $A$ is solid if and only if $\ell_\infty \subset (A,A)$; 2. $A \subset (B,H)$ if and only if $B \subset (A,H)$;
3. the solid core $s(A)$ of $A$ is the largest solid space contained in $A$, and moreover $s(A)=(\ell_\infty,A)$; 4. the solid hull $S(A)$ of $A$ is the smallest solid space containing $A$;  5. If $X$ is solid, $(A,X)=(S(A),X)$ and $(X,A)=(X,s(A))$.

The results on \cite{BT1} and \cite{BT2} contain a characterization of the solid hull of
$H_v^\infty$, if $v$ satisfies a condition $(b)$, see \eqref{bbb} below.
This general  characterization is in
terms of a numerical sequence $(m_n)_{n=0}^\infty$, which depends on the weight and was
studied by the second named author in \cite{L2}. However, given a concrete weight on the disc like
\bea
v(r) = \exp \big( - a / (1-r)^b \big)
\eea
the calculation of the numbers $m_n$ is not an easy matter, and in \cite{BT2}, this
was only done in the case $ 0 < b \leq 2$. In Theorem \ref{th4.1} we calculate these numbers
and thus determine the solid hull for weights
\bea
v(r) = w(r) \exp \big( - a / (1-r)^b \big)
\eea
with any $a,b>0$, where $w$ is a differentiable positive function with some growth
restriction. Moreover, the same theorem also contains the analogous characterization
of the solid core. In Section \ref{secNEW} we show how these results can be used
also in the case of the variant
\bea
v(r) = \exp \big( - a / (1-r^2)^b \big)
\eea
of the weight.  The general characterization of solid cores for weights satisfying
condition $(b)$ is given in Theorem \ref{th2.2}. As explained above, further interesting related results are presented in Section \ref{sec5}.

Bennet, Stegenga and Timoney in their paper \cite{BST} determined the solid hull and the solid core of the weighted spaces $H^\infty_v(\bbD)$  in the case the weight $v$ is doubling. Exponential weights $v(r)=\exp(-a/(1-r)^b)$ with $a,b>0$ are not doubling.
Not much seems to be known about multipliers and solid hulls of weighted spaces of analytic functions on the unit disc in the case of exponential weights. Hadamard multipliers of certain weighted space $H^1_a(\alpha), \alpha >0,$ were completely described by Dostani\'c in \cite{Dos} (see also Chapter 13 in \cite{JVA}). Other aspects of weighted spaces of analytic functions on the unit disc with exponential weights, like integration operators or Bergman projections, have been investigated recently by Constantin, Dostani\'c, Pau, Pavlovi\'c, Pel\'aez and R\"atty\"a, among others; see \cite{CP}, \cite{Dos2}, \cite{PP}, \cite{Pav} and \cite{PR}. The solid hull and multipliers on spaces of analytic functions on the disc has been investigated by many authors. In addition to  \cite{BST}, we mention for example  \cite{AS}, \cite{BG}, \cite{BP}, \cite{JP}, the books \cite{JVA} and \cite{Pav-book} and the many references therein.

Spaces of type $H_v^\infty(\bbC)$ and $H^\infty_v(\bbD)$ appear
in the study of growth conditions of
analytic functions and have been investigated in various articles since
the work of Shields and Williams, see {\it e.g.} \cite{BBG},\cite{BBT},
\cite{L1}, \cite{L2}, \cite{SW} and the references therein.

In the case of a "standard" weight $v_\alpha (z) )= (1-|z|^2)^\alpha$, where
$\alpha \geq 0$, we denote for every
$H_\alpha^\infty:= H^\infty_\alpha (\bbD) := H_{v_\alpha}^\infty$.
The solid hull $S\big( H^\infty_\alpha \big)$ of  $H^\infty_\alpha $ is
known: it equals
$$ 
S\big( H^\infty_\alpha \big) = \Big\{ (b_m)_{m=0}^\infty \, : \,
\sup_{n \in \bbN_0} \Big( \sum_{m=2^n}^{2^{n+1}-1} |b_m|^2 (m+1)^{-2\alpha}
 \Big)^{1/2} < \infty \  \Big\} .
$$ 
This is Theorem 8.2.1 of \cite{JVA}.
Moreover, the solid core $s\big( H^\infty_\alpha \big)$ can also be characterized,
see Theorem 8.3.4 of \cite{JVA}:
$$ 
s\big( H^\infty_\alpha \big)
 = \Big\{ (b_m)_{m=0}^\infty \, : \,
\sup_{n \in \bbN_0} \Big( \sum_{m=2^n}^{2^{n+1}-1} |b_m| (m+1)^{-\alpha}
 \Big) < \infty \  \Big\} . 
$$ 

\section{Solid hull and core for weights with condition $(b)$.}
\label{sec2}

In this section we consider the quite large class of weights on $\bbD$ or $\bbC$ 
satisfying the regularity condition $(b)$. For such a weight, the solid hull was 
found in the paper \cite{BT2}, and we determine the solid core here.

\begin{definition}\label{condition(b)}{\rm
Let $r_n \in ]0,R[ $ be a global maximum point of the function
$r^m v(r)$ for any $m > 0$. The weight $v$ satisfies the \em{condition}
$(b)$ if there exist numbers $b> 2$, $K > b$ and $ 0 < m_1 < m_2 < \ldots$ with
$\lim_{n \rightarrow \infty} m_n = \infty$ such that
\bea
b \leq \left( \frac{r_{m_n}}{r_{m_{n+1}}} \right)^{m_n}
\frac{v(r_{m_n})}{v(r_{m_{n+1}})}, \left( \frac{r_{m_{n+1}}}{r_{m_{n}}} \right)^{m_{n+1}}
\frac{v(r_{m_{n+1}})}{v(r_{m_{n}})} \leq K.     \label{bbb}
\eea}
\end{definition}

\begin{remark}
\label{remLB}{\rm
(1) The second named author introduced the following \textit{condition (B)} on the weight $v$ in \cite{L2}:
$$
\forall b_1>0 \ \exists b_2 >1 \ \exists c>0 \ \forall m,n :
$$
$$
\left(\frac{r_m}{r_n}\right)^m \frac{v(r_m)}{v(r_n)} \leq b_1 \ \ {\rm and} \ \ |m-n| \geq c \Rightarrow \left(\frac{r_n}{r_m}\right)^n \frac{v(r_n)}{v(r_m)} \leq b_2.
$$
It was observed in Remark 2.7 of \cite{BT1} that if a weight $v$ satisfies condition (B), then it also satisfies condition $(b)$ for some $b> 2$, $K > b$ and $ 0 < m_1 < m_2 < \ldots$.

(2) As a consequence of this observation and Section 2 in \cite{L2}, the following weights satisfy condition $(b)$: For $R=1$,

\noindent $(i)$ $v(r)=(1-r)^{\alpha}$ with $\alpha >0$, which are the standard weights on the disc, and

\noindent $(ii)$ $v(r)= \exp(-(1-r)^{-1})$.

\noindent More examples can be seen in Example \ref{ex1.5}. For $R=+\infty$,

\noindent $(i)$ $v(r)= \exp(-r^p)$ with  $p > 0$,

\noindent $(ii)$ $v(r)= \exp( -\exp r)$,  and

\noindent $(iii)$ $v(r)= \exp\big(- (\log^+ r  )^p\big) $, where $p \geq 2$ and $\log^+ r
= \max(\log r,0)  $.}

\end{remark}

We recall the result \cite{BT2}, Theorem 2.1 for $\bbD$, or  
\cite{BT1}, Theorem 2.5 for entire functions:

\BET
\label{th2.1}
If the weight $v$ satisfies $(b)$, we have
\bea \hspace{3cm} S(H_v^\infty) = \bigg\{ (b_m)_{m=0}^\infty : \sup_n
v(r_{m_{n}}) \Big( \sum_{m_n < m \leq m_{n+1} } |b_m|^2 r_{m_n}^{2m} \Big)^{1/2}
< \infty \bigg\} . \label{2.4}
\eea
\ENT

Now let us prove the following statement.

\BET
\label{th2.2}
For a weight $v$ satisfying $(b)$, we have
\bea \hspace{3cm} s(H_v^\infty) = \bigg\{  (b_m)_{m=0}^\infty : \sup_n
v(r_{m_{n}}) \Big( \sum_{m_n < m \leq m_{n+1} } |b_m| r_{m_n}^{m} \Big) < \infty \bigg\} .
\label{2.6}
\eea
\ENT

Proof.
For a holomorphic function $f$ with $f(z) = \sum_{n=0}^{\infty} a_n z^n$, we let $h_f$ denote the function defined by
\[ h_f(z)=\sum_{n=0}^{\infty} |a_n| z^n \ \ \mbox{ for all } z. \]
It is easy to see that  the  solid core $s(H_v^\infty)$ of $H_v^\infty$ coincides with the set
\[ \{ f : f \mbox{ holomorphic}, h_f \in H_v^\infty \}. \]

Now, let $f \in H_v^\infty$, $f(z) = \sum_{n=0}^{\infty} b_n z^n$. If
\bea
h_f (z) = \sum_{n=0}^{\infty} |b_n| z^n , \label{2.7}
\eea
belongs to $H_v^\infty $, it is easily seen that
\[\Vert h_f \Vert_v = \sup_{0 < r < R} v(r) \sum_{n=0}^{\infty} |b_n|r^n.\]
This implies
\bea
\label{(2.1)}  \sup_n  v(r_{m_{n}})
\Big( \sum_{m_n < m \leq m_{n+1} } |b_m| r_{m_n}^m \Big) \leq \Vert h_f\Vert_v.
\eea
Thus the solid core is contained in the right-hand side of (\ref{2.6}). Now we proceed with the reverse inclusion.

According to \cite{L2}, Proposition 5.2., there are numbers $\beta_m \in [0,1]$ and a constant $c>0$ (independent of $f$) such that
\[ \Vert  h_f\Vert_v \leq c\sup_n\sup_{r_{m_{n-1}} \leq |z| \leq r_{m_{n+1}}} v(z) \Big|\sum_{m_{n-1} < m \leq m_{n+1}  }\beta_m |b_m| z^m \Big|.   \]
This implies
\bea
\label{(2.2)}
 \Vert h_f\Vert_v \leq c\sup_n\sup_{r_{m_{n-1}} \leq r \leq r_{m_{n+1}}} v(r) \Big( \sum_{m_{n-1} < m \leq m_{n+1}  } |b_m| r^m \Big).
\eea
For $r_{m_{n-1}} \leq r \leq r_{m_{n+1}}$ we have
\begin{eqnarray}
 & & v(r) \sum_{m_{n-1} < m \leq m_{n}  } |b_m| r^m \
\roweq
\frac{	v(r)}{v(r_{m_{n-1}})} v(r_{m_{n-1}})\sum_{m_{n-1} < m \leq m_{n}  } |b_m| r^m_{m_{n-1}}
\Big(\frac{r}{r_{m_{n-1}}} \Big)^m
\rowleq
\Big(\frac{r}{r_{m_{n-1}}}  \Big)^{m_n}   \frac{	v(r)}{v(r_{m_{n-1}})} v(r_{m_{n-1}})\sum_{m_{n-1} < m \leq m_{n}  } |b_m| r^m_{m_{n-1}}
\rowleq
\Big(\frac{r_{m_n}}{r_{m_{n-1}}}  \Big)^{m_n}   \frac{	v(r_{m_n})}{v(r_{m_{n-1}})} v(r_{m_{n-1}})\sum_{m_{n-1} < m \leq m_{n}  } |b_m| r^m_{m_{n-1}}
\rowleq
K v(r_{m_{n-1}}) \sum_{m_{n-1} < m \leq m_{n}  } |b_m| r^m_{m_{n-1}}.
\label{(2.3)}
\end{eqnarray}                                                                                       Here we  used that $r_{m_n}$ is a global maximum point for $r^{m_n}v(r)$.

Similarly, for     $r_{m_{n-1}} \leq r \leq r_{m_{n}}$   we have
\begin{eqnarray}
& & v(r) \sum_{m_{n} < m \leq m_{n+1}  } |b_m| r^m
\roweq
\frac{	v(r)}{v(r_{m_{n}})} v(r_{m_{n}})\sum_{m_{n} < m \leq m_{n+1}  } |b_m| r^m_{m_{n}}
\Big(\frac{r}{r_{m_{n}}} \Big)^m
\rowleq
\Big(\frac{r}{r_{m_{n}}}  \Big)^{m_n}
\frac{	v(r)}{v(r_{m_{n}})} v(r_{m_{n}})\sum_{m_{n} < m \leq m_{n+1}  } |b_m| r^m_{m_{n}}
\rowleq	
v(r_{m_{n}})\sum_{m_{n} < m \leq m_{n+1}  } |b_m| r^m_{m_{n}}.
\label{(2.4)}
\end{eqnarray}
Finally, if  $r_{m_{n}} \leq r \leq r_{m_{n+1}}$  then
\begin{eqnarray}
& & v(r) \sum_{m_{n} < m \leq m_{n+1}  } |b_m| r^m
\roweq
\frac{	v(r)}{v(r_{m_{n}})} v(r_{m_{n}})\sum_{m_{n} < m \leq m_{n+1}  } |b_m| r^m_{m_{n}}
\Big(\frac{r}{r_{m_{n}}} \Big)^m
\rowleq
\Big(\frac{r}{r_{m_{n}}}  \Big)^{m_{n+1}}   \frac{	v(r)}{v(r_{m_{n}})} v(r_{m_{n}})\sum_{m_{n} < m \leq m_{n+1}  } |b_m| r^m_{m_{n}}
\rowleq
\Big(\frac{r_{m_{n+1}}}{r_{m_{n}}}  \Big)^{m_{n+1}}
\frac{	v(r_{m_{m+1}})}{v(r_{m_{n}})} v(r_{m_{n}})\sum_{m_{n} < m \leq m_{n+1}  } |b_m| r^m_{m_{n}}
\rowleq
K v(r_{m_{n}})\sum_{m_{n} < m \leq m_{n+1}  } |b_m| r^m_{m_{n}}.
\label{(2.5)}
\end{eqnarray}
Hence, according to \eqref{(2.2)}, with \eqref{(2.3)}, \eqref{(2.4)} and \eqref{(2.5)},
\[
\Vert h_f\Vert_v \leq c 2K \sup_n  v(r_{m_{n}})  \Big(
\sum_{m_n < m \leq m_{n+1} } |b_m| r_{m_n}^m \Big).
\]
This together with (2.1) yields \eqref{2.6}. $ \Box$

\section{A weight of the form $w(r) \exp(-a/(1-r)^{b})$.}

In this section we only deal with weights defined on the unit disc $\bbD$.
Our purpose is to improve the results of \cite{BT2} by calculating the solid
hulls and cores for a larger class of concrete examples, namely, for weights of the form
\bea
\label{(3.1.1)}
v(z) = w(r) \exp\Big( - \frac{a}{(1-r)^b} \Big), \ \ z \in \mathbb{D},
\eea
where  $a,b > 0$ are  given constants and  $w: [0,1[ \to ]0,\infty[$ is a differentiable  
function, extended to $\bbD$ by $w(z) = w(|z|)$. We remark
the examples  in \cite{BT2} only contain the case $b \leq 2$, $w \equiv 1$.

We will prove
\BET
\label{th4.1}
Let $w'(r)/w(r)$ be a decreasing function and assume that  there are $n_0 > 0$ and
$\alpha \in ]0, 1+b/2[$ such that
\bea	\label{(3.1.2)}    & &    (1-r)^{\alpha} \frac{w'(r)}{w(r)} \mbox{ is bounded on }
[0,1[  \\
\label{(3.1.3)}  & &     \frac{1}{e} \leq \frac{w(1- (\frac{a}{bn^2})^{1/b})}{w(1- (\frac{a}{b(n+1)^2})^{1/b})} \leq e \mbox{ for } n \geq n_0.
\eea
Then, the solid hull of $H_v^\infty$ is equal to
\beas
\bigg\{(b_m)_{m=0}^\infty: \sup_n
w \Big(1- \Big( \frac{a}{bn^2}\Big)^{1/b} \Big)e^{-bn^2}
\bigg(  \!\!\!\!\! \sum_{{m \in \mathbb{N}} \atop {m_n < m \leq m_{n+1}}}  \!\!\!\!\!
|b_m|^2 \Big( 1- \Big( \frac{a}{bn^2}\Big)^{1/b}\Big)^{2m} \bigg)^{1/2} < \infty \bigg\} ,
\eeas
where
\beas m_n = b \left( \frac{b}{a} \right)^{1/b} n^{2+2/b} -bn^2-
\Big(1-\left(\frac{a}{bn^2} \right)^{1/b}\Big)
\frac{w'(1-\left(\frac{a}{bn^2} \right)^{1/b})}{w(1-\left(\frac{a}{bn^2} \right)^{1/b})}
\eeas
(which reduces to
$ m_n = b^{1+1/b} a^{-1/b}  n^{2+2/b} -bn^2$, if  $w \equiv 1$).

Moreover, the solid core of $H_v^\infty$ is equal to
\beas
\bigg\{ (b_m)_{m=0}^\infty:
\sup_n w\Big(1- \Big( \frac{a}{bn^2} \Big)^{1/b} \Big)  e^{-bn^2}
\!\!\!\!\!  \!\!\! \sum_{{m \in \mathbb{N}} \atop {m_n < m \leq m_{n+1}}}   \!\!\!\!\!
|b_m|\Big(1-  \Big( \frac{a}{bn^2}\Big)^{1/b}\Big)^{m}  < \infty \bigg\}.
\eeas
\ENT

We postpone the proof a bit and consider some remarks and examples.
Let us 
start with the following  quite trivial observation which, however,  is
very useful to simplify the presentations of the solid hulls and cores.

\BEL
\label{lem1.4}
Let $1  \leq p < \infty$, and let $(K_n)_{n \in \bbN_0}$,
 $(\tilde K_n)_{n \in \bbN_0}$, $(L_{m})_{n \in \bbN_0}$
and $(\tilde L_{m})_{n \in \bbN_0}$ be sequences of positive numbers.
Assume that there are given two increasing, unbounded  sequences  $(m_n)_{n \in \bbN_0}$ and
$(\tilde m_n)_{n \in \bbN_0}$ of positive real numbers
such that
\bea
m_n < \tilde m_n < m_{n+1} \ \ \ \forall \ n \in \bbN_0
\eea
and such that for some constants $C > c > 0$, for all $n \in \bbN$,
\bea
& & c K_n^p L_m \leq \tilde K_{n-1}^p \tilde L_m \leq
C K_n^p L_m  \ \ \ \ \forall \, m \, \mbox{with} \ m_n < m \leq \tilde m_n
\nonumber \\
& & c K_n^p L_m \leq \tilde K_{n}^p \tilde L_m \leq
C K_n^p L_m  \ \ \ \ \forall \, m \, \mbox{with} \ \tilde m_n < m \leq  m_{n+1}.
\label{1.40}
\eea
%
Then, we have
\bea
& & \Big\{ (b_m)_{m=0}^\infty \, : \,
\sup_{n \in \bbN_0} K_n \Big( \sum_{m_n < m \leq m_{n+1}} |b_m|^p L_m
\Big)^{1/p} < \infty \  \Big\}
\roweq
\Big\{ (b_m)_{m=0}^\infty \, : \,
\sup_{n \in \bbN_0} \tilde K_n  \Big( \sum_{\tilde m_n < m \leq \tilde m_{n+1}}
|b_m|^p \tilde L_m  \Big)^{1/p} < \infty \  \Big\}
\eea
\ENL

Proof. We have, by \eqref{1.40},
\bea
& &
\sup_{n \in \bbN_0}  \Big( \sum_{m_n < m \leq m_{n+1}}  K_n^p |b_m|^p L_{m}
\Big)^{1/p}
\rowleq
\sup_{n \in \bbN_0} \bigg(   \Big( \sum_{m_n < m \leq \tilde m_n}
 K_n^p |b_m|^p L_{m}  \Big)^{1/p}
+   \Big( \sum_{\tilde m_n < m \leq m_{n+1}}
 K_n^p |b_m|^p L_{m}  \Big)^{1/p} \bigg)
\rowleq
\frac1c \sup_{n \in \bbN_0} \bigg(   \Big( \sum_{m_n < m \leq \tilde m_n}
\tilde K_{n-1} ^p|b_m|^p \tilde L_{m}  \Big)^{1/p}
+   \Big( \sum_{\tilde m_n < m \leq m_{n+1}}
\tilde K_n^p |b_m|^p \tilde L_{m}  \Big)^{1/p} \bigg)
\rowleq
\frac2c \sup_{n \in \bbN_0}   \Big( \sum_{\tilde m_n < m \leq \tilde m_{n+1}}
\tilde K_n^p |b_m|^p \tilde L_{m}  \Big)^{1/p} . \nonumber
\eea
The converse inequality can be shown in the same way, using the other inequalities in
\eqref{1.40}.  \ \ $\Box$

\bigskip

It is obvious that the representations of the solid hull and core of a weighted space
are by no means unique: the sequence $m_n$ and even the coefficients and exponents
can be chosen in many ways. We discuss this in the first example.

\BEE
\label{ex1.5}
$(i)$: $a=b=1$, $w=1$. Here $m_n = n^4-n^2$. However,  in \cite{BT2}
the representation of this solid hull was found with the more simple numbers
$\tilde m_n = n^4$ instead: 

\beas
& & S(H_v^\infty) = \Big\{ (b_m)_{m=0}^\infty : \sup_n e^{-n^2}
\Big(  \sum_{m = n^4 +1 }^{(n+1)^4} 
|b_m|^2 \big(1- n^{-2}\big)^{2m} \Big)^{1/2} < \infty \Big\} ,
\\
&&
s(H_v^\infty) = \Big\{  (b_m)_{m=0}^\infty : \sup_n e^{-n^2}
\sum_{m = n^4 +1 }^{(n+1)^4}
|b_m|  \big(1- n^{-2}\big)^{m}
 < \infty \Big\} . \eeas
Let us verify that the condition \eqref{1.40} of
Lemma \ref{lem1.4} is satisfied with $p=2$; the proof for the case $p=1$ follows by taking square roots.
We may obviously assume $n \geq 2$. For $m$ with $m_n < m \leq \tilde m_n$, i.e.,
\bea
n^4 - n^2 < m \leq n^4,      \label{3.60}
\eea
we have $K_n = e^{-n^2} = \tilde K_n$, $L_m = (1- n^{-2})^{2m}$ and $\tilde L_m = (1- (n-1)^{-2})^{2m}$,
hence,
\bea
K_n^2 L_m & = & e^{-2 n^2} \Big( 1 - \frac1{n^2} \Big)^{2m} \geq
e^{-2 n^2} \Big( 1 - \frac1{n^2} \Big)^{n^2 2n^2  }
\rowgeq
C  e^{-2 n^2}  e^{-2 n^2} = C e^{-4 n^2}   ,   \label{3.62} \\
K_n^2 L_m & \leq  & 
e^{-2 n^2} \Big( 1 - \frac1{n^2} \Big)^{2n^4 - 2n^2 }
\rowleq 
e^{-2 n^2}  \Big( 1 - \frac1{n^2} \Big)^{n^2 2n^2 } \Big( 1 - \frac1{n^2} \Big)^{ - 2n^2 }
\rowleq C  e^{-4 n^2} \Big( 1 - \frac1{n^2} \Big)^{- 2 n^2  }
\leq
C'  e^{- 4 n^2}    ,   \label{3.64}
\eea
Furthermore, we write 
$$
n^4 = (n-1)^4 + 4 (n-1)^3+ \rho(n) 
= (n-1)^2(n-1)^2 + 4 (n-1)^3+ \rho(n),
$$  
where  $\rho(n)= 6n^2 - 8n +3$. Using the trivial estimates  $\rho(n)- n^2 \geq - 60 (n-1)^2$
and   $\rho(n)\leq 50 (n-1)^2$ for all $n \geq 2$, we obtain 
\bea
\tilde K_{n-1}^2 \tilde L_m & = & e^{-2(n-1)^2} \Big( 1 - \frac1{(n-1)^2} \Big)^{2 m}
\leq e^{-2(n-1)^2} \Big( 1 - \frac1{(n-1)^2} \Big)^{2n^4 - 2 n^2 }
\roweq
e^{-2(n-1)^2} \Big( 1 - \frac1{(n-1)^2} \Big)^{2(n-1)^2(n-1)^2 + 8 (n-1)^3 
+ 2 ( \rho(n) -  n^2) }
\rowleq
e^{-2(n-1)^2} \bigg( \Big( 1 - \frac1{(n-1)^2} \Big)^{(n-1)^2} \bigg)^{2(n-1)^2 + 8 (n-1)
- 120 }
\rowleq
C e^{- 4 (n-1)^2}  e^{-8(n-1)} =
C e^{- 4 n^2 + 8 n -4}  e^{-8(n-1)} \leq C' e^{- 4 n^2 } .   \label{3.66}
\eea
Similarly,
\bea
\tilde K_{n-1}^2 \tilde L_m & = & e^{-2(n-1)^2} \Big( 1 - \frac1{(n-1)^2} \Big)^{2m  }
\geq e^{-2(n-1)^2} \Big( 1 - \frac1{(n-1)^2} \Big)^{2n^4 }
\roweq
e^{-2(n-1)^2} \Big( 1 - \frac1{(n-1)^2} \Big)^{2(n-1)^2(n-1)^2 + 8 (n-1)^3 + 2 \rho(n)  }
\rowgeq
e^{-2(n-1)^2} \bigg( \Big( 1 - \frac1{(n-1)^2} \Big)^{(n-1)^2} \bigg)^{2(n-1)^2 + 8 (n-1)
+100 }
\rowgeq
C e^{- 4 (n-1)^2}  e^{-8(n-1)} 
\geq C' e^{- 4 n^2 } .   \label{3.66a}
\eea
We thus see that the first pair of inequalities \eqref{1.40} holds. The second one is
trivial since for  $\tilde m_n < m \leq  m_{n+1}$, i.e.,
\beas
n^4 < m \leq (n+1)^4 - (n+1)^2,     
\eeas
we have $L_m = (1- n^{-2})^{2m} = \tilde L_m$.

\bigskip

\noindent$(ii)$: $a=1$, $b=2$, $w(r) = 1-r$. Here
\beas m_n = 2^{3/2} n^3- 2n^2 + 2^{1/2}n -1 \eeas
and	
\beas
& & S(H_v^\infty) = \Big\{ (b_m)_{m=0}^\infty : \sup_n \frac{e^{-2n^2}}{\sqrt{2}n} \Big(
\!\!\!\!\! \sum_{{m \in \mathbb{N}} \atop {m_n < m \leq m_{n+1}}} \!\!\!\!\!
|b_m|^2 \big(1- (\sqrt{2}n)^{-1}\big)^{2m}
\Big)^{1/2} < \infty \Big\} ,
\\ & &
s(H_v^\infty) = \Big\{(b_m)_{m=0}^\infty : \sup_n \frac{e^{-2n^2}}{\sqrt{2}n}
\!\!\!\!\! \sum_{{m \in \mathbb{N}} \atop {m_n < m \leq m_{n+1}}} \!\!\!\!\!
|b_m|  \big(1-(\sqrt{2} n)^{-1}\big)^{m}
< \infty \Big\} . \eeas

\bigskip

\noindent$(iii)$: $a=b=1$, $w(r) = (1- \log(1-r))^{-1}$. Here, a direct calculation
yields
\bea
m_n = n^4-n^2 + \frac{n^2-1}{1+ \log(n^2)} ,   \label{3.tx}
\eea
but we can again use Lemma \ref{lem1.4} with $\tilde m_n = n^4$, 
$K_n = e^{-n^2}\big( 1+ \log(n^2) \big)^{-1} = \tilde K_n$, $L_m = (1- n^{-2})^{2m}$ and $\tilde L_m = (1- (n-1)^{-2})^{2m}$, since the calculation \eqref{3.62}--\eqref{3.66a}
shows  that both the expressions $K_n^2 L_m$ and $\tilde K_{n-1}^2 L_m$ 
are proportional to 
\bea
\frac{e^{-4 n^2}}{\big( 1+ \log(n^2) \big)^2} \label{3.tz}
\eea 
for all $m$ with $ m_n < m \leq \tilde m_n$.  (To see this, we observe that in comparison 
with  \eqref{3.62}--\eqref{3.66a}, $\tilde K_n$ only has the new factor $\big( 1+ \log(n^2) 
\big)^{-1} =: g_n$ for which $g_n$ and $g_{n-1}$ are proportional, and that $m_n$ of 
\eqref{3.tx} satisfies $n^4 - n^2 < m_n < n^4$ so that $m$ in \eqref{3.tz} falls into
the interval considered also in  \eqref{3.62}--\eqref{3.66a}). 
Moreover, of course $K_n^2 L_m = \tilde K_{n}^2 L_m$
for all $m$ with $ \tilde m_n < m \leq  m_{n+1}$. Thus, we have
\beas
& & S(H_v^\infty) = \Big\{ (b_m)_{m=0}^\infty : \sup_n  \frac{e^{-n^2}}{1+ \log(n^2)}
\Big( \sum_{m = n^4+1}^{(n+1)^4} 
|b_m|^2 \big(1- n^{-2}\big)^{2m}
\Big)^{1/2} < \infty \Big\} , \\
& &
s(H_v^\infty) = \Big\{ (b_m)_{m=0}^\infty : \sup_n  \frac{e^{-n^2}}{1+ \log(n^2)}
 \sum_{m = n^4+1}^{(n+1)^4} 
|b_m|  \big(1- n^{-2}\big)^{m}
 < \infty \Big\} .
\eeas

\noindent$(iv)$:  $a=b=1$, $w(r)= \exp(- \log^2(1-r)) $. Here we have $w'(r)/w(r) = 2 (1-r)^{-1}
\log(1-r)$. It is easily seen that \eqref{(3.1.2)}  and \eqref{(3.1.3)} are satisfied.
We obtain  $m_n= n^4 - n^2 + 4(n^2-1) \log(n)$ and
\beas
& & S(H_v^\infty) =
\Big\{ (b_m)_{m=0}^\infty: \sup_n \exp(-4 \log^2(n)-n^2)
\Big(
\!\!\!\!\! \sum_{{m \in \mathbb{N}} \atop {m_n < m \leq m_{n+1}}} \!\!\!\!\!
|b_m|^2 \big(1- n^{-2}\big)^{2m}
\Big)^{1/2} < \infty \Big\} , \\
& &
s(H_v^\infty) =
\Big\{ (b_m)_{m=0}^\infty: \sup_n \exp(-4 \log^2(n)-n^2)
\!\!\!\!\! \sum_{{m \in \mathbb{N}} \atop {m_n < m \leq m_{n+1}}} \!\!\!\!\!
|b_m|  \big(1- n^{-2}\big)^{m}
< \infty \Big\} . \eeas
\ENE

\BER
Fix  $m > 1$ and  put
\beas
f(r) = r^m v(r) = r^m w(r) \exp \Big(- \frac{a}{(1-r)^b} \Big).
\eeas
Due to the continuity of $f$ and the fact that $f(0) = f(1) = 0 \leq f(r)$, $ r \in ]0,1[$, the
function $f$ has a global maximum on $]0,1[$. It  is easily seen that $r \in ]0,1[$ is a zero of
$f'$ if and only if
\bea
\label{(3.3.1)}  m = ab \frac{r}{(1-r)^{b+1}} - r \frac{w'(r)}{w(r)}.
\eea
Since $ -r w'(r)/w(r)$ is assumed to be increasing, the right-hand side of \eqref{(3.3.1)}
is strictly
increasing in $r$. Hence \eqref{(3.3.1)} has exactly one solution, denoted by $r_m$, which is
the unique global maximum of $f$. In particular, if
\bea
M = ab m^{1+1/b}\Big( 1- \Big( \frac{1}{m} \Big)^{1/b} \Big) -
\Big( 1- \Big( \frac{1}{m} \Big)^{1/b} \Big)
\frac{w'\Big( 1- \Big( \frac{1}{m} \Big)^{1/b} \Big)}{
w\Big(1- \Big( \frac{1}{m} \Big)^{1/b} \Big)}  \label{NNM}
\eea
for some $m > 1$,  then  
\bea
r_M =1 -\Big( \frac{1}{m} \Big)^{1/b}.  \label{MMM}
\eea
\ENR

Proof of Theorem \ref{th4.1}. 
$1^\circ$. We first consider that case  $b \geq 1$.

$a)$ Some estimates.
If $1 \leq x \leq y$ and $ 0 < \beta < 1 $ then the mean value theorem yields
\bea
\label{(3.4.1)} y^\beta-x^{\beta} \leq \beta
x^{\beta-1}(y-x) \ \ \mbox{ and } \ \ x^{\beta} - y^{\beta} \leq \beta y^{\beta-1}(x-y).
\eea
Moreover we use
\bea
\label{(3.4.2)}  1+x \leq e^x \ \ \ \mbox{ for all } x \in \mathbb{ R }.
\eea
Now let $1 < m \leq k$ and  define $M$ as in \eqref{NNM} and $K$ in the same
way with $k$ replacing $m$. 
Then
$$ 
r_M = 1 -\frac{1}{m^{1/b}} \ \mbox{ and } \ \ r_K = 1 -\frac{1}{k^{1/b}}
$$
and we can rewrite \eqref{NNM} as 
\bea
M= ab m^{1+1/b} r_M -  r_M  \frac{w'(r_M)}{w(r_M)} \ , \ \ 
K= ab k^{1+1/b} r_K -  r_K  \frac{w'(r_K)}{w(r_K)} .   \label{NNN}
\eea
Then, with \eqref{(3.4.1)},
 \eqref{(3.4.2)} we obtain
\begin{eqnarray}
 \frac{r_K}{r_M} & =  & 
\frac{1 - \left( \frac{1}{k} \right)^{1/b}}{1 - \left( \frac{1}{m} \right)^{1/b}} =
 1 + \frac{\left( \frac{1}{m} \right)^{1/b} - \left( \frac{1}{k} \right)^{1/b}}{1 -\left( \frac{1}{m} \right)^{1/b}} \nonumber  \\
 & \leq   &
\exp \bigg( \frac{1}{b} \Big( \frac{1}{k} \Big)^{1/b - 1}
\Big( \frac1m - \frac1k \Big) \frac{1}{r_M}  \bigg)=
  \exp\bigg( \frac{1}{b} \  \frac{k-m}{k^{1/b} m } \ \frac1{r_M} \bigg)
  \label{(3.4.3)} 
\end{eqnarray}  
 and
\begin{eqnarray}
 \frac{r_M}{r_K} & =  & \frac{1 - \left( \frac{1}{m} \right)^{1/b}}{1 - \left( \frac{1}{k} \right)^{1/b}}
= 
 1 - \frac{\left( \frac{1}{m} \right)^{1/b} - \left( \frac{1}{k} \right)^{1/b}}{1 -\left( \frac{1}{k} \right)^{1/b}} \nonumber  \\
& \leq   &
\exp\bigg( -\frac{1}{b} \ \frac{k-m}{m^{1+ 1/b}} \ 
\frac{1}{r_K} \bigg). \label{(3.4.4)}
\end{eqnarray}

Now we write
$$
\Big( \frac{r_K}{r_M} \Big)^K \frac{v(r_K)}{v(r_M)}
= \exp \bigg( K \log \Big( \frac{r_K}{r_M} \Big) \bigg) \frac{w(r_K)}{w(r_M)} \exp(- a(k-m) )
$$
which  in view of \eqref{NNN}, \eqref{(3.4.3)} is bounded by
$$
\frac{w(r_K)}{w(r_M)} \ \exp
\bigg(a \frac{k(k-m)}{m} \ \frac{r_K}{r_M}
- \frac{1}{b} \ \frac{r_K}{r_M} \ \frac{w'(r_K)}{w(r_M)} \
\frac{k-m}{k^{1/b}m} - a(k-m)  \bigg)
$$
$$
= 
\frac{w(r_K)}{w(r_M)} \ \exp
\bigg( - a(k-m) + a  \bigg( \frac{k- k^{1 - 1/b}}{m-m^{1-1/b}}\bigg)(k-m)
+ c_1(k,m)
$$
\bea
= \frac{w(r_K)}{w(r_M)} \  \exp \bigg( a \frac{(k-m)^2}{m - m^{1- 1/b}}
\bigg(1 - \frac{k^{1-1/b} - m^{1-1/b}}{k-m}
\bigg) + c_1(k,m) \bigg)
\label{(3.4.7)}
\eea
where
\bea
\label{(3.4.8)}
c_1(k,m) = - \frac{1}{b} \ \frac{r_K}{r_M} \  \frac{w'(r_K)}{w(r_K)} (1-r_K)
\frac{k-m}{m}.
\eea
Using \eqref{(3.4.4)} instead of \eqref{(3.4.3)} we get with the help of
\eqref{(3.4.1)} 
$$
 \Big( \frac{r_M}{r_K} \Big)^K \frac{v(r_M)}{v(r_K)}
$$
$$
\leq
\frac{w(r_M)}{w(r_K)}\exp\bigg( a(k-m) - a (k-m) \left( \frac{k}{m} \right)^{1/b} 
+ c_2(k,m) \bigg)
$$
$$
= \frac{w(r_M)}{w(r_K)}   \exp\left(a(k-m) \Big(\frac{m^{1/b}-k^{1/b}}{m^{1/b}} \Big) + c_2(k,m) \right)
$$
\bea
\leq
\frac{w(r_M)}{w(r_K)}\exp\left(-\frac{a}{b}(k-m)^2\frac{1}{k^{1-1/b}m^{1/b}}+c_2
(k,m)\right)  \label{(3.4.9)}
\end{eqnarray}
  where
\bea
\label{(3.4.10)}    c_2(k,m) =	 \frac{1}{b} \left(\frac{1}{m}\right)^{1/b} 
\frac{w'(r_K)}{w(r_K)} \ \frac{k-m}{k}  .
\eea
Similarly, \eqref{(3.4.3)} implies
$$
 \Big( \frac{r_K}{r_M} \Big)^M \frac{v(r_K)}{v(r_M)}
$$
$$
\leq \frac{w(r_K)}{w(r_M)} \exp\left(-a(k-m) + a \Big( \frac{m}{k} \Big)^{1/b}(k-m) + c_3(k,m) \right)
$$
$$
= \frac{w(r_K)}{w(r_M)}   \exp\left(-a(k-m) \Big(\frac{k^{1/b}-m^{1/b}}{k^{1/b}} \Big)
+ c_3(k,m) \right)
$$
\bea
\leq
\frac{w(r_K)}{w(r_M)} \exp\left( -\frac{a}{b} (k-m)^2 \frac{1}{k} + c_3(k,m) \right)
\label{(3.4.11)}
\end{eqnarray}
with
\bea
\label{(3.4.12)}    c_3(k,m) = - \frac{1}{b} \  \frac{w'(r_M)}{w(r_M)}
\ \frac{k-m}{m} \ (1-r_K).
\eea
Finally, \eqref{(3.4.4)} implies
$$
\Big( \frac{r_M}{r_K} \Big)^M \frac{v(r_M)}{v(r_K)} 
$$
$$ 
=  \frac{w(r_M)}{w(r_K)} \exp\left(a(k-m) -a \frac{m}{k}(k-m) 
\frac{ r_M}{r_K}  + c_4(k,m) \right)
$$
\bea
= \frac{w(r_M)}{w(r_K)}\exp\left(a\frac{(k-m)^2}{k-k^{1-1/b}}
\Big( 1 - \frac{k^{1-1/b}-m^{1-1/b}}{k-m} \Big) + c_4(k,m) \right)
\label{(3.4.13)}
\end{eqnarray}
with
\bea
\label{(3.4.14)}    c_4(k,m) = \frac{1}{b} \frac{r_M}{r_K} \ 
\frac{w'(r_M)}{w(r_M)} \ \frac{k-m}{m^{1/b} k} .
\eea

$b)$ The parameters $m_n$. Now put for every $n \in \bbN$ large enough
\bea
\label{(3.4.15)} j_n =  \frac{b}{a} n^2 
\eea
and in the above calculations choose
$$
m = j_{n} \ \ , \ \ k =j_{n+1} .
$$
We denote the numbers in \eqref{NNN} by 
$$
m_n = M \ \ , \ \ m_{n+1} = K
$$
so that the following relations  hold, by \eqref{MMM}:
$$ 
r_{m_n} = 1-\frac{1}{j_n^{1/b}}=1- \left(\frac{a}{bn^2}\right)^{1/b} \ , \ \ 
r_{m_{n+1}} = 1- \left(\frac{a}{b(n+1)^2}\right)^{1/b}.
$$ 

$c)$ Final estimates. 
With \eqref{(3.4.8)} we obtain
$$ 
c_1(j_{n+1},j_n) = - \frac{1}{b} \Big( \frac{r_{m_{n+1}}}{r_{m_n}} \Big) \frac{w'(r_{m_{n+1}})}{w(r_{m_{n+1}})}(1-r_{m_{n+1}}) \ \frac{2n+1}{n^2} 
$$ 
By assumption \eqref{(3.1.2)} there is $d > 0$ and $0 < \alpha < 1+b/2$ with
$$
\frac{w'(r_{m_{n+1}})}{w(r_{m_{n+1}})} \leq \frac{d}{(1-r_{m_{n+1}})^{\alpha}}.
$$ 
Hence
$$
|c_1(j_{n+1},j_n)|
\leq  \frac{d}{b}\Big( \frac{r_{m_{n+1}}}{r_{m_n}} \Big)(1-r_{m_{n+1}})^{1-\alpha}
\Big(\frac{2n+1}{n^2} \Big) 
$$
$$ 
= \frac{d}{b}\Big(\frac{r_{m_{n+1}}}{r_{m_n}} \Big) \Big(\frac{a}{b}\Big)^{(1-\alpha)/b}(n+1)^{2(\alpha-1)/b}\Big(\frac{2n+1}{n^2} \Big).
$$ 
Since $ 2(\alpha -1)/b < 1$ we obtain	
$$ 
\lim_{n \rightarrow \infty}c_1(j_{n+1},j_n) = 0
$$ 
By assumption \eqref{(3.1.3)} we have $  w(r_{m_{n+1}})/w(r_{m_{n}})   \leq e$. So, using \eqref{(3.4.7)} and \eqref{(3.4.15)} we see that there is a constant $K_0$ with
\bea
\label{(3.4.21)}
\Big(  \frac{r_{m_{n+1}}}{r_{m_n}} \Big)^{m_{n+1}} \frac{v(r_{m_{n+1}})}{v(r_{m_{n}})} \leq K_0 \ \ \ \mbox{ for all } n.
\eea
This follows from the fact that, for $k= j_{n+1}$ and $m = j_n$, the expression
in \eqref{(3.4.7)},
$$
\frac{(k-m)^2}{m -m^{1/b}} \ \frac{1- k^{1-1/b}-m^{1-1/b}}{k-m} ,
$$
remains uniformly bounded for all $n$.

To obtain a lower estimate of 
$$ 
\left(  \frac{r_{m_{n+1}}}{r_{m_n}} \right)^{m_{n+1}} \frac{v(r_{m_{n+1}})}{v(r_{m_{n}})}
$$ 
consider \eqref{(3.4.9)} and \eqref{(3.4.10)}. Exactly as before we see that
\bea
\label{(3.4.22)}  \lim_{n \rightarrow \infty}c_2(j_{n+1},j_n) = 0
\eea
For $k=j_{n+1}$ and $m=j_n$ we obtain
$$
\Big(\frac{a}{b}\Big) \frac{(k-m)^2}{k^{1-1/b}m^{1/b}} = \frac{(2n+1)^2}{(n+1)^{2-2/b}n^{2/b}}
$$ 
which tends to $4$ as $n \rightarrow \infty$. Together with \eqref{(3.4.22)}
we find $n_0$ such that
$$
- \Big(\frac{a}{b}\Big)\frac{(j_{n+1}-j_n)^2}{j_{n+1}^{1-1/b}j_n^{1/b}} + c_2(j_{n+1},j_n)   \leq -2
\ \ \ \mbox{ for } n \geq n_0.
$$
Since by assumption $w(r_{m_n})/w(r_{m_{n+1}}) \leq e$ the estimate \eqref{(3.4.9)} implies
$$
\Big(  \frac{r_{m_{n}}}{r_{m_{n+1}}} \Big)^{m_{n+1}} \frac{v(r_{m_{n}})}{v(r_{m_{n+1}})} \leq \frac{1}{e},
$$
hence
\bea
\label{(3.4.23)}
2 < e \leq  \Big(  \frac{r_{m_{n+1}}}{r_{m_n}} \Big)^{m_{n+1}} \frac{v(r_{m_{n+1}})}{v(r_{m_{n}})}  \ \ \ \ \mbox{ for } n \geq n_0.
\eea
Repeating the preceding arguments using \eqref{(3.4.11)}, \eqref{(3.4.12)},
\eqref{(3.4.13)}, \eqref{(3.4.14)} instead of \eqref{(3.4.7)}, \eqref{(3.4.8)},
\eqref{(3.4.9)}, \eqref{(3.4.10)} we see that
$$
\lim_{n \rightarrow \infty}c_3(j_{n+1},j_n) = \lim_{n \rightarrow \infty}c_4(j_{n+1},j_n)
= 0
$$
and there are $n_1$, $N_1$ with
\bea
\label{(3.4.24)}
 2< e \leq  \Big(  \frac{r_{m_{n}}}{r_{m_{n+1}}} \Big)^{m_{n}} \frac{v(r_{m_{n}})}{v(r_{m_{n+1}})} \leq K_1 \ \ \ \mbox{ for } n \geq n_1.
\eea
Then the assertion of the theorem in the case $b \geq 1$ follows from 
\eqref{(3.4.21)}, \eqref{(3.4.23)},
\eqref{(3.4.24)} and \cite{BT2},  Theorem 2.1.

$2^\circ$. We prove Theorem \ref{th4.1} in the case $0 < b < 1$.
Here we use, for $ \gamma > 1$, $ 0 \leq x \leq y$,
$$
y^\gamma-x^{\gamma} \leq \gamma
y^{\gamma-1}(y-x) \ \ \mbox{ and } \ \ x^{\gamma} - y^{\gamma} \leq \gamma x^{\gamma-1}(x-y).
$$
We obtain, instead of \eqref{(3.4.3)} and  \eqref{(3.4.3)},
\begin{eqnarray*}
\frac{1 - \left( \frac{1}{k} \right)^{1/b}}{1 - \Big( \frac{1}{m} \Big)^{1/b}}
& =  &
1 + \frac{\left( \frac{1}{m} \right)^{1/b} - \left( \frac{1}{k} \right)^{1/b}}{1 -\left(
\frac{1}{m} \right)^{1/b}}
\\     & \leq   &
\exp\bigg( \frac{1}{b}\left( \frac{1}{m} \right)^{1/b} \frac{k-m}{\Big(1 -\left( \frac{1}{m} \right)^{1/b}\Big)k} \bigg)
\end{eqnarray*}
and
\begin{eqnarray*}
\frac{1 - \left( \frac{1}{m} \right)^{1/b}}{1 - \left( \frac{1}{k} \right)^{1/b}}
& =   &
1 - \frac{\left( \frac{1}{m} \right)^{1/b} - \left( \frac{1}{k} \right)^{1/b}}{1 -\left( \frac{1}{k} \right)^{1/b}}
\\ & \leq   &
\exp\bigg( -\frac{1}{b}\Big( \frac{1}{k} \Big)^{1/b} \frac{k-m}{\Big(1 -\left( \frac{1}{k}
\right)^{1/b}\Big)m} \bigg).
\end{eqnarray*}
Then the theorem follows by repeating the same arguments as in the preceding section.
\ \ $\Box$

\section{A weight of the form  $v_2(z)=\exp(-a/(1-r^2)^b)$.}
\label{secNEW}
As a consequence of the preceding discussion we consider here the weight
\[ v_2(z)=\exp\Big( \frac{-a}{(1-r^2)^b}\Big) \]
for given constants $a,b > 0$. We compare $v_2$ with the weight $v_1(z)= \exp(-a/(1-r)^b)$ of the preceding section (with $w \equiv 1$).

Put 
\[ 
A = \Big\{ f \in H_{v_2}^{\infty} : f(z) = \sum_{k=0}^{\infty} a_{2k}z^{2k} \ \mbox{ for some } a_{2k} \Big\} 
\]
and
\[ 
B = z \cdot A = \Big\{ g \in H_{v_2}^{\infty} : g(z) = \sum_{k=0}^{\infty} a_{2k+1}z^{2k+1} \ \mbox{ for some } a_{2k+1} \Big\}. 
\]
Moreover, let ${T_1},{T_2} : H_{v_1}^{\infty} \rightarrow H_{v_2}^{\infty}$ be the maps with
\[ 
({T_1}h)(z) = h(z^2) \ \ \ \mbox{ and } \ \ \ ({T_2}h)(z)= z h(z^2), \ \ h\in H_{v_1}^{\infty}, z \in \mathbb{D}. 
\]

\BEP
\label{prop41}
The operator 
${T_1}$ maps $H_{v_1}^{\infty}$ isometrically onto $A$. The map ${T_2}$ is a contractive operator 
from $H_{v_1}^{\infty}$ 	onto $B$. Moreover, we have  	\[ H_{v_2}^{\infty} = A \oplus B. 
\]
\ENP

Proof. The map ${T_1}$ is certainly an isometry into $A$ and ${T_2}$ is a contractive operator into 
$B$. To show the surjectivity, let $f  \in B$, say 
\bea
f(z) = \sum_{k=0}^{\infty} a_{2k+1}z^{2k+1}. \label{(4.1)}
\eea
Then put $ h(z) = \sum_{k=0}^{\infty} a_{2k+1}z^{k}$. In view of \eqref{(4.1)}, since the series representing $f(z)/z$ converges uniformly on compact subsets  of $\mathbb{D}$, we have
\[ c:= \sup_{|z| \leq 1/2} \frac{|f(z)|}{|z|}  < \infty.  \]
Hence,
$$
\Vert h \Vert_{v_1} =  \sup_{|z| < 1} |h(z)|v_1(z) = \sup_{|z| <1} |h(z^2)|v_1(z^2) 
=  \sup_{|z| < 1} \frac{|f(z)|}{|z|} v_2(z) 
$$
\bea 
= \max \bigg( \sup_{|z| \leq 1/2} \frac{|f(z)|}{|z|} v_2(z), \sup_{1/2<|z| < 1} 
\frac{|f(z)|}{|z|} v_2(z) \bigg) 
\eea
We obtain $h \in H_{v_1}^{\infty}$ and clearly ${T_2}h=f$. This shows that ${T_2}$ maps $H_{v_1}^{\infty}$ onto $B$. Similarly we see that ${T_1}$ maps $H_{v_1}^{\infty}$ onto $A$.

Now consider the operator $P$ with $(Pf)(z) = (f(z)+f(-z))/2$ for $f \in H_{v_2}^{\infty}$ and $z \in \mathbb{D}$.
$P$ is a contractive projection from $H_{v_2}^{\infty}$ onto $A$. We clearly get 
$( id-P)(H_{v_2}^{\infty}) = B$.
Hence $H_{v_2}^{\infty} = A \oplus B$. $\Box$

\bigskip

Now we take the numbers $m_n$ of Theorem \ref{th4.1} for $ w \equiv 1$, i.e.
\[ m_n = \frac{{b}^{1+1/b}}{a^{1/b}} \ n^{2+2/b}-bn^2.\]
Let $[s]$ denote the largest integer which is smaller than or equal to $ s$.

\BET
\label{th42}
	The solid hull of $H_{v_2}^{\infty}$ is equal to 
\[	\Big\{ (b_m) : \sup_n e^{-bn^2} \Big(
\!\!\!\!\! \sum_{{m \in \mathbb{N}} \atop {2[m_n]+1 < m \leq 2[m_{n+1}]+1}} \!\!\!\!\!
|b_m|^2 \Big(1- \Big( \frac{a}{bn^2}\Big)^{1/b}\Big)^{2[m/2]} \ \ \ \ \ \Big)^{1/2} < \infty \Big\}.\]	
Moreover, the solid core of $H_{v_2}^{\infty}$ is equal to
\[	\Big\{ (b_m) : \sup_n e^{-bn^2} 
\!\!\!\!\! \sum_{{m \in \mathbb{N}} \atop {2[m_n]+1 < m \leq 2[m_{n+1}]+1}} \!\!\!\!\!
 |b_m| \Big(1- \Big( \frac{a}{bn^2}\Big)^{1/b}\Big)^{[m/2]} < \infty \Big\}.\]	
\ENT

Proof. Using Proposition \ref{prop41} and Theorem \ref{th4.1} we see that $(b_m) \in 
S(H_{v_2}^{\infty})$ if and only if
\bea
\sup_n e^{-bn^2} \Big(
\!\!\!\!\! \sum_{{m \in \mathbb{N}} \atop {m_n < m \leq m_{n+1}}} \!\!\!\!\!
|b_{2m}|^2 \Big(1- \Big( \frac{a}{bn^2}\Big)^{1/b}\Big)^{2m}\Big)^{1/2} < \infty 
\label{(4.2)}
\eea
and
\bea
\sup_n e^{-bn^2} \Big(
\!\!\!\!\! \sum_{{m \in \mathbb{N}} \atop {m_n < m \leq m_{n+1}}} \!\!\!\!\!
|b_{2m+1}|^2 \Big(1- \Big( \frac{a}{bn^2}\Big)^{1/b}\Big)^{2m}\Big)^{1/2} < \infty. 
\label{(4.3)}
\eea
If $m \in \mathbb{N}$ and $m_n < m \leq m_{n+1}$ then $[m_n]+1 \leq m \leq [m_{n+1}]$.  
We obtain
$2[m_n]+2 \leq 2m \leq 2 [m_{n+1}] $ and $2[m_n]+3 \leq 2m+1 \leq 2[m_{n+1}]+1$.
Hence \eqref{(4.2)} and \eqref{(4.3)} are equivalent to
\[ \sup_n e^{-bn^2} \Big(
\!\!\!\!\! \sum_{{m \in \mathbb{N}} \atop {2[m_n]+1 < m \leq 2[m_{n+1}]+1}} \!\!\!\!\!
|b_m|^2 \Big(1- \Big( \frac{a}{bn^2}\Big)^{1/b}\Big)^{2[m/2]} \ \ \ \ \ \Big)^{1/2} < \infty. \] 
The proof for the solid core is the same. \ \ $\Box$

\section{Maximal solid cores and minimal solid hulls.}
\label{sec5}


In this section we prove some general results on the
relations of Schauder bases and solid hulls and cores
for $H_v^\infty$-spaces. We refer the reader to \cite{LT} for terminology about bases in Banach spaces. If the monomials form a Schauder basic sequence in  $H_v^\infty$, then obviously the condition of being solid is related with the property of $\{ z^m \}$ being an
unconditional basis. Another, related fact will be proven in
Theorem \ref{th5.2}. We also show the unexpected fact that for some special
weights, $H_v^\infty$ is solid.

\BEP
\label{prop5.1}
We have
$S(H_v^\infty) = H_v^\infty$ if and only if $s(H_v^\infty)= H_v^\infty$.
\ENP

Proof. If 	$S(H_v^\infty) = H_v^\infty$ then $h_f \in H_v^\infty$ (see \eqref{2.7})
for all $f \in H_v^\infty$. Hence   $s(H_v^\infty) = H_v^\infty$.

\bigskip

Now assume $s(H_v^\infty) = H_v^\infty$ and take $g \in S(H_v^\infty)$ with $g(z) = \sum_{k=0}^{\infty} b_k z^k$. There is
 $f \in H_v^\infty$ with  $f(z) = \sum_{k=0}^{\infty} a_k z^k$ and $|b_k| \leq |a_k|$ for all $k$. Since by assumption $h_f \in H_v^\infty$
  we obtain
  \[\Vert g\Vert_v \leq \sup_{r} v(r)\sum_{k=0}^{\infty} |a_k| r^k = \Vert h_f\Vert_v < \infty \]
  which implies $g \in H_v^\infty$. Hence 	$S(H_v^\infty) = H_v^\infty$. $\Box$

\bigskip

{\bf Example.} Consider the weight $v(r) = \exp(- \log^2(r))$ on the complex plane $\mathbb{C}$.
According to \cite{L1}, Theorem 2.5., there is a constant $ d > 0$ such that for every  $f \in
H_v^\infty$ with   $f(z) = \sum_{k=0}^{\infty} a_k z^k$ we have
\[
\sup_k (|a_k| \exp(k^2/4)) \leq \Vert f\Vert_v  \leq d  \sup_k (|a_k| \exp(k^2/4)).
\]
Then clearly $h_f \in H_v^\infty$. Indeed, let $(h_f)_n$ be the partial sums of $h_f$, i.e.
$(h_f)_n(z)=  \sum_{k=0}^n |a_k|z^k$. Then $(h_f)_n \in H_v^\infty$, $(h_f)_n \rightarrow h_f$
pointwise on $\mathbb{C}$   and
\[
\Vert h_f\Vert_v \leq \sup_n \Vert(h_f)_n\Vert_v \leq d
 \sup_k (|a_k| \exp(k^2/4)) \leq d \Vert f \Vert_v < \infty.
\]
Hence $S(H_v^\infty) = H_v^\infty= s(H_v^\infty)$.

\bigskip

Recall that we denote by  $H_{v}^0$  the closure of the polynomials in $H_v^\infty$. We put $\Lambda = \{z^k : k =0,1,2, \ldots \}$.

\BET
\label{th5.2}
If $S(H_v^\infty)= H_v^\infty$ then $\Lambda$ is a Schauder basis of $H_{v}^0$.
\ENT

We prove Theorem \ref{th5.2} at the end of this section. At first we state
\BEC
\label{cor5.3}
In the case of analytic functions on the disc $\mathbb{D}$, one always has $S(H_v^\infty(\mathbb{D})) \not= H_v^\infty(\mathbb{D})$ and $s(H_v^\infty(\mathbb{D})) \not= H_v^\infty(\mathbb{D})$.
\ENC

{\bf Proof.} According to \cite{L1}, Theorem 2.2., $\Lambda$ is never a basis for $H_{v}^0(\bbD)$. This proves Corollary \ref{cor5.3}  in view of Theorem \ref{th5.2} $\Box$

\bigskip

For the proof of Theorem \ref{th5.2} we need two lemmas.

\BEL
\label{lem5.4}
 (i)	Fix $m \in \mathbb{N}$ and $\epsilon > 0$. Then there is $r_0 < R$ such that, for every $f$ with
$f(z) = \sum_{k=0}^{m} a_k z^k$, we have
\[
\sup_{r_0 \leq |z| < R} |f(z)| v(z) \leq \epsilon \Vert f \Vert_v.
\]
(ii) Fix $0 \leq r_1 < R$ and $\epsilon > 0$. Then there is $n \in \mathbb{N}$ such that, for any $g \in H_v^\infty$
 with $g(z) = \sum_{k=n}^{\infty} a_k z^k$, we have
 \[      \sup_{0 \leq |z| \leq r_1} |g(z)| v(z) \leq \epsilon \Vert g\Vert_v. \]
\ENL

{\bf Proof.} (i) Fix $r < R$. Then we clearly have
\[ |a_k| \leq \frac{\Vert f\Vert_v}{r^kv(r)} \ \ \mbox{ for all } k= 0,1, \ldots, m. \]
We obtain
 \[ |f(z)| v(z) \leq \sum_{k=0}^m \Big( \frac{|z|}{r} \Big)^k \frac{v(|z|)}{v(r)} \Vert f\Vert_v. \]
Find $r_0> 0$ such that
\[
\Big( \frac{|z|}{r} \Big)^k  \frac{v(|z|)}{v(r)} \leq  \frac{\epsilon}{m+1}
\]
whenever $|z| > r_0$. This is possible since, by assumption, $ \lim_{|z| \rightarrow R} |z|^k v(|z|) =0$
for all $k$. This implies (i).
\vspace{4mm} \\
(ii)   Fix $r > r_1$ and consider $g(z) = \sum_{k=n}^{\infty}a_kz^k$. We have $|a_k| \leq \Vert g \Vert _v / (r^k v(r)) $ for all $k$. This implies, if $|z| \leq r_1$,
\[
|g(z)|v(z) \leq \sum_{k=n}^{\infty} \frac{|z|^kv(|z|)}{r^kv(r)} \Vert g\Vert_v \leq \sum_{k=n}^{\infty} \Big(\frac{r_1}{r} \Big)^k \frac{v(0)}{v(r)} \Vert g\Vert_v .
\]
We find $n$ so large that
\[
\sum_{k=n}^{\infty} \Big(\frac{r_1}{r} \Big)^k \frac{v(0)}{v(r)} \leq \epsilon
\]
which proves the lemma. $\Box$

\BEL
\label{lem5.5} 
Let $f = \sum_{j=0}^\infty a_j z^j$ be an  analytic function on the disc,
let $m_1 < m_2 < \ldots$ be indices and $f_n(z) =
\sum_{j= m_n+1}^{m_{n+1}} a_j z^j$. 
Then there is a subsequence $(f_{n_k})_{k=0}^\infty $ such that
\bea 
\sup\limits_{k \in \bbN} \Vert f_{n_k}\Vert_v \leq  2\Vert  \sum_{k=0}^\infty 
f_{n_k}\Vert_v.  \label{5.76}
\eea
\ENL

We remark that for any  subsequence $(f_{n_k})_{k=0}^\infty $, the sum 
$\sum_k f_{n_k}$ on the right-hand side of \eqref{5.76} 
is the Taylor series of an analytic function on the disc, so the sum
converges at least uniformly on compact subsets of $\bbD$; if the sum does not belong to $H_v^\infty$, its norm
is infinity and the inequality \eqref{5.76} becomes a triviality. 

\bigskip

{\bf Proof.} Use Lemma \ref{lem5.4} and induction to find a subsequence $(f_{n_k})$ and radii $r_1 < r_2 < \ldots$
such that
$ |f_{n_k}(z)|v(|z|) \leq 3^{-k} \Vert 	f_{n_k}\Vert_v$ whenever $|z| \leq r_k$ or $ |z| \geq r_{k+1}$.
Hence
\bea
\label{5.1}
\Vert 	f_{n_k}\Vert_v = \sup_{r_k \leq |z| \leq r_{k+1}} |f_{n_k}(z)|v(|z|).
\eea
By the remark above, we may assume that  $\sum_k f_{n_k} \in H_v^\infty$.
Fix $j$.  If $ r_j \leq |z| \leq r_{j+1}$ we obtain
	\begin{eqnarray*}
\Vert  \sum_kf_{n_k}\Vert_v & \geq & |\sum_kf_{n_k}(z)|v(z) \\
	                 & \geq & |f_{n_j}(z)|v(z) - \sum_{k \not= j}|f_{n_k}(z)|v(z) \\
	                 & \geq & |f_{n_j}(z)|v(z) - \sum_{k \not= j} \frac{1}{3^k}\Vert f_{n_k}\Vert_v \\
	                 & \geq & |f_{n_j}(z)|v(z) - \frac{1}{2} \sup_k \Vert f_{n_k}\Vert_v.
	\end{eqnarray*}
In view of \eqref{5.1}
this implies
\[
\Vert \sum_kf_{n_k}\Vert_v  \geq  \Vert f_{n_j}\Vert_v - \frac{1}{2} \sup_k \Vert f_{n_k}\Vert_v  \ \ \mbox{ for all } j
\]
and hence
\[
\Vert \sum_kf_{n_k}\Vert_v  \geq   \frac{1}{2} \sup_k \Vert f_{n_k}\Vert_v
\]
which proves the lemma. $\Box$
	
\bigskip
	
{\bf Proof of Theorem \ref{th5.2}}. For any subset  $N$ of $\mathbb{N}$, let $T_N$ be the
operator with	$T_N(\sum_{k =0}^{\infty}a_kz^k) = \sum_{k \in {N}}a_kz^k $. If $S(H_v^\infty)
= H_v^\infty = s(H_v^\infty)$ then 	$T_N(H_v^\infty) \subset H_v^\infty$. The closed graph theorem implies that $T_N$ is bounded.
	
Now   let $P_n$ be the Dirichlet projections, i.e. $P_n(\sum_{k =0}^{\infty}a_kz^k) =
\sum_{k=0}^na_kz^k$. Assume that $\Lambda$ is not a basis for $H_{v}^0$. Then
the $P_n$ are not 	uniformly bounded. By the uniform boundedness theorem we obtain a function
$f \in H_{v}^0$ such 	that $\sup_n\Vert P_n(f)\Vert_v= \infty$. Hence we can find a
subsequence $P_{n_m}$ with 	$\lim_{m \rightarrow \infty}\Vert (P_{n_{m+1}} -P_{n_m})
(f)\Vert_v = \infty$. Put $f_m=(P_{n_{m+1}} -P_{n_m})(f)$. Then, $ \sum_m f_{m} \in 
H_v^\infty$, since this sum is of the form $T_Nf$ for some subset $N $ of $\bbN$. 
We apply Lemma \ref{lem5.5} to find a subsequence $f_{m_k}$ such that
\[
\sup_k \Vert f_{m_k} \Vert_v \leq 2 \Vert  \sum_k f_{m_k}\Vert_v.
\]
The left hand side of this inequality is infinite while the function on the right-hand side is 
again of the form 	$T_{\tilde N} f$ for some $\tilde N \subset \bbN$ and thus has finite norm as an element 
$H_v^\infty$. 
So we arrive at a contradiction. 	Therefore $\Lambda$ is a basis of $H_{v}^0$. $\Box$
	
\bigskip

\noindent \textbf{Acknowledgements.} The research of Bonet was partially
supported by the project  MTM2016-76647-P. The research of Taskinen was
partially supported by the V\"ais\"al\"a Foundation of the Finnish Academy
of Sciences and Letters.

\vspace{.5cm}

\noindent \textbf{Authors' addresses:}%
\vspace{\baselineskip}%

Jos\'e Bonet: Instituto Universitario de Matem\'{a}tica Pura y Aplicada IUMPA,
Universitat Polit\`{e}cnica de Val\`{e}ncia,  E-46071 Valencia, Spain

email: jbonet@mat.upv.es \\

Wolfgang Lusky: FB 17 Mathematik und Informatik, Universit\"at Paderborn, D-33098 Paderborn, Germany.

email: lusky@uni-paderborn.de \\

Jari Taskinen: Department of Mathematics and Statistics, P.O. Box 68,
University of Helsinki, 00014 Helsinki, Finland.

email: jari.taskinen@helsinki.fi


\begin{thebibliography}{99}

\bibitem{AS} J.M. Anderson, A.L. Shields, Coefficient multipliers of Bloch functions, Trans. Amer. Math. Soc. 224 (1976), 255-265.

\bibitem{BST} G. Bennet, D.A. Stegenga, R.M. Timoney, Coefficients of Bloch and Lipschitz functions, Illinois J. Math. 25 (1981), 520-531.

\bibitem{BBG} K.D.\ Bierstedt, J.\ Bonet, A.\ Galbis, Weighted spaces
	of holomorphic functions on bounded domains, Michigan
Math.\ J.\ 40 (1993), 271--297.

\bibitem{BBT} K.D.\ Bierstedt, J.\ Bonet, J.\ Taskinen, Associated
	weights and spaces of holomorphic functions, Studia Math.\
	127 (1998), 137--168.

\bibitem{BG} O. Blasco, A. Galbis, On Taylor coefficients of entire functions integrable against exponential weights, Math. Nachr. 223 (2001), 5-21.

\bibitem{BP} O. Blasco, M. Pavlovic, Coefficient multipliers on Banach spaces of analytic functions, Rev. Mat. Iberoam. 27 (2011),  415-447.


\bibitem{BT1} J.Bonet, J.Taskinen, Solid hulls of weighted Banach spaces of entire functions, to appear in Rev. Mat. Iberoamericana. arXiv: 1607.02237v1

\bibitem{BT2} J.Bonet, J.Taskinen, Solid hulls of weighted Banach spaces of analytic functions on the unit disc with exponential weights, to appear in Ann.Acan.Sci.Fenn. arXiv:1702.00145v1.

\bibitem{CP} O. Constantin, J.A. Pel\'aez, Boundedness of the Bergman projection on $L_p$-spaces with exponential weights, Bull. Sci. Math. 139 (2015), no. 3, 245-268.

\bibitem{Dos} M.R. Dostani\'c, Multipliers in the space of analytic functions with exponential mean growth, Asymptot. Anal. 65 (2009), no. 3-4, 191--201.

\bibitem{Dos2} M-R. Dostani\'c, Integration operators on Bergman spaces with exponential weight, Rev. Mat. Iberoam. 23 (2007), no. 2, 421-436.

\bibitem{JP} M. Jevti\'c, M.  Pavlovi\'c, On the solid hull of the Hardy-Lorentz space,
Publ. Inst. Math. (Beograd) (N.S.) 85(99) (2009), 55-61.

\bibitem{JVA} M. Jevti\'c, D. Vukoti\'c, M. Arsenovi\'c, Taylor Coefficients and Coefficient Multipliers of Hardy and Bergman-Type Spaces, RSME Springer Series, Volume 2. Springer 2016.


\bibitem{LT} J.Lindenstrauss, L.Tzafriri, Classical Banach spaces I, Springer, Berlin, 1977.

\bibitem{L1} W.Lusky,  On the Fourier series of unbounded harmonic functions, J. of the Lond. Math. Soc. (2) 61, 568-580, (2000)

\bibitem{L2} W.Lusky, On the isomorphism classes of weighted spaces of harmonic and holomorphic functions, Studia Math. 175, 19-45, (2006)

\bibitem{PP} J. Pau, J.A. Pel\'aez, Volterra type operators on Bergman spaces with exponential weights. Topics in complex analysis and operator theory, 239-252, Contemp. Math., 561, Amer. Math. Soc., Providence, RI, 2012.

\bibitem{Pav} M. Pavlovi\'c, On harmonic conjugates with exponential mean growth,
Czech. Math. J. 49 (1999), 733--742.

\bibitem{Pav-book} M. Pavlovi\'c, Function classes on the unit disc.
An introduction, De Gruyter Studies in Mathematics, 52. De Gruyter, Berlin, 2014. xiv+449 pp.

\bibitem{PR} J.A. Pel\'aez, J. R\"atty\"a, Weighted Bergman spaces induced by rapidly increasing weights,
Mem. Amer. Math. Soc. 227 (2014), no. 1066, vi+124 pp.

\bibitem{SW} A.L. Shields, D.L. Williams, Bounded projections, duality
	and multipliers in spaces of analytic functions, Trans.
Amer. Math. Soc. 162 (1971), 287-302.

\end{thebibliography}
\end{document}